\documentstyle[12pt]{article}

  \newcommand{\const}{\rm const}
  \newcommand{\Var}{\rm Var}
  
  \newcommand{\Law}{\rm Law}

  \newcommand{\Sub}{\rm  Sub}
  \newcommand{\diam}{\rm diam}
  \newcommand{\Cov}{\rm  Cov}

\begin{document}

   \begin{center}

 \ \  {\bf Improved  Monte-Carlo method for solving of  } \\

\vspace{4mm}

 \ \ {\bf integral Fredholm's equations of a second kind, } \\

 \vspace{4mm}

  \ \ {\bf   with confidence regions in the uniform norm} \\

\vspace{4mm}

  \ \ {\bf   Ostrovsky E., Sirota L.}\\

\vspace{4mm}

 Israel,  Bar-Ilan University, department of Mathematic and Statistics, 59200, \\

\vspace{4mm}

e-mails: \  eugostrovsky@list.ru \\
sirota3@bezeqint.net \\

\vspace{4mm}

 \ \ {\bf Abstract} \\

\vspace{4mm}

 \end{center}

 \ \  We offer in this article some modification of Monte-Carlo method for solving of a linear integral Fredholm's
equation of a second kind (Fredholm's well posed problem). \par

 \ We prove that the rate of convergence of offered method is optimal under natural conditions still in an uniform norm, and
construct an asymptotic as well as non-asymptotic confidence region, again in the uniform norm. \par

\vspace{4mm}

 \ \ {\it  Key words and phrases:} Kernel, Linear integral Fredholm's equation of a second kind, Monte-Carlo
method, random variables, natural distance, Central Limit Theorem in the space of continuous functions, contraction, Kroneker's and
ordinary degree of integral operator, Gaussian random field, covariation function, tail of distribution,
metric entropy, entropic integral, ordinary and subgaussian norm and space, norm of linear operator, compact operators, Neuman
series, uniform norm, spectral radius,  asymptotic and non-asymptotic confidence region, random variable and random vector,
variance, Dependent Trial Method (DTM).\par

\vspace{4mm}

 \ Mathematics Subject Classification 2000. Primary 42Bxx, 4202, 68-01, 62-G05,
90-B99, 68Q01, 68R01; Secondary 28A78, 42B08, 68Q15.

\vspace{4mm}

 \section{\ Definitions. Notations. Previous results. Statement of problem.}

 \vspace{4mm}

 \ Let $ \  T = \{t\}  \  $ be compact  metrisable space equipped with probabilistic Borelian non-trivial complete measure: $ \  \mu(T) = 1. \ $
The completeness  of the measure implies that each non-empty open set $ \ A \ $ has a positive measure: $ \  \mu(A) > 0. \ $\par

 \ We consider a linear Fredholm's type equation in the space $ \ C(T) \ $  of continuous functions equipped as ordinary with the uniform norm

$$
\ ||g|| := \max_{t \in T} |g(t)|
$$
of the form

$$
z(t) = f(t) + \int_T K(t,s) \ z(s) \ \mu(ds) = f(t) + \int K(t,s) \ z(s) \ \mu(ds) =
$$

$$
f(t) + K[z](t), \ t,s \in T, \eqno(1.1)
$$
where the  kernel  $ \  K = K(t.s)  \  $ is continuous function of two  variables $ \ (t,s) \in T^2, \ $ so that
$ \  K[\cdot] \  $ is a linear integral (compact) operator of the form

$$
 K[z](t) \stackrel{def}{=} \int_T K(t,s) \ z(s) \ \mu(ds),  \ z(\cdot) \in C(T),
$$
 and  with continuous {\it non-zero} ``free'' member $ \ f = f(t), \ t \in T. \ $ \par

 \ {\it We will write for brevity in the sequel}

$$
\int g(t) \ \mu(dt)  \stackrel{def}{=}  \int_T g(t) \ \mu(dt), \ g(\cdot) \in C(T)
$$
or at last for the function $ \ g (\cdot) \in L_1(T, \mu). \ $ \par

 \  Put

$$
\rho = \rho[K] \stackrel{def}{=} \max_{t \in T} \int_{T} |K(t,s)| \ \mu(ds), \eqno(1.2)
$$

$$
\overline{\rho} = \overline{\rho}[K] \stackrel{def}{=}   \max_{t,s \in T} |K(t,s)|. \eqno(1.2a)
$$

 \ Obviously, $  0 \le \rho[K] \le \overline{\rho}[K]. \ $ \par

 \ The value $ \ \rho = \rho[K] \ $ is nothing more that the norm of operator $ \ K \ $ acting in the space $ \ C(T): \ $

$$
\rho[K]  = ||K||(C(T) \to C(T)) = \sup_{0 \ne g \in C(T)} \left\{ \frac{||K[g]||}{||g||} \right\}. \eqno(1.3)
$$

\vspace{4mm}

 \ {\it  We suppose  hereafter that the operator   } $ \ K[\cdot] \ $ {\it satisfies in addition the contraction principle:}

$$
\rho[K]  = ||K||(C(T) \to C(T)) := \sup_{0 \ne g \in C(T)} \left\{ \frac{||K[g]||}{||g||} \right\}  < 1. \eqno(1.4)
$$

\vspace{4mm}

 \ The classical expression for solution $ \ z(t) \ $ may be written  by virtue of the uniformly convergent iterations of fixed point method (Neuman series):
 $ \ z(t) = \lim_{d \to \infty} z_d(t), \ $  where $ \  z_0(t) = f(t), \ $ and for $ \ d = 0,1,2,\ldots $

 $$
 z_{d+1}(t) = f(t) +\int_T K(t,s) \ z_d(s) \ \mu(ds) = f(t) + K[z_d](t). \eqno(1.5)
 $$

 \ The iteration $ \ z_d(t) \ $ has a form $ \ z_0(t) = K^0[f](t) = f(t), \ $

$$
 z_1(t) = K[f](t) =\int_T K(t,s) \ f(s) \ \mu(ds),
$$

 and in the case when $  d \ge 2 \  $

$$
  z_d(t) = K^d[f](t) = \int_{T^d} K_d(t, \vec{s}) \  \prod_{l=1}^d \mu(ds_l)  =
$$

$$
\int_{T^d} K_d(t, s_1, s_2, \ldots, s_d) \ \cdot f(s_d)  \cdot  \prod_{l=1}^d \mu(ds_l)  =
$$

$$
\int_T \int_T  \ldots \int_T K(t,s_1) \cdot  \left[\prod_{r=1}^{d - 1}  K \left(s_r, s_{r+1} \right) \right] \cdot f(s_d) \cdot
\left[ \prod_{l=1}^d \mu(ds_l) \right].  \eqno(1.6)
$$

 \ Here $ \  K^d[\cdot] \  $ denotes the $ \ d^{th} \ $ power of the (integral) operator $ \ K[\cdot], \ $ in contradiction to the ordinary power of
the function $ \ K^2(t,s), \ $

$$
\vec{s} = \vec{s}_d = (s_1,s_2, \ldots, s_d), \   K_d(t, \vec{s}) = K_d(t, s_1, s_2, \ldots, s_d) :=
$$

$$
K(t, s_1) \cdot  \left[\prod_{r=1}^{d - 1}  K \left(s_r, s_{r+1} \right) \right].
$$

 \ The truncated sum

$$
 z^{(M)}(t)  \stackrel{def}{=} f(t) + \sum_{d=1}^{M} z_d(t), \ M = 2,3,\ldots
$$
gives a following error estimate

$$
||z(\cdot) - z^{(M)}(\cdot) || \le ||f|| \cdot \frac{\rho^{M+1}}{1 - \rho}.
$$

 \ In order to calculate each integral $ \ z_d(t), \ $ we offer the classical Monte-Carlo method, which is named  as "Dependent Trial Method", see
[9], [10], [19], [20].  Indeed, let us  introduce a $ \ d \ - \ $ dimensional random  vector $ \  \vec{\xi}^{(d)}:  \  $

$$
\vec{\xi}^{(d)} = \left\{ \xi^d_1,  \xi^d_2, \ldots, \xi^d_d \right\},
$$
where  all the random variables $ \  \xi^d_k \ $ are (common) independent and have the distributions  $ \ \mu: \ $

$$
{\bf P} \left( \ \xi^d_k \in A \  \right) = \mu(A).
$$

 \ Let also $ \  \vec{\xi}^{(d)}_i \  $ be independent copies of the r.v.  $ \ \vec{\xi}^{(d)}.  \  $ The consistent as $ \ n(d) \to \infty, \ $
 here and in the sequel $ \ n(d), \ d = 1,2,\ldots \  $ are "Great"  integer numbers, unbiased
Monte-Carlo approximation $ \  x_{d, n(d)}(t) \  $ for the integral $ \ x_d(t) \ $  may has a form

$$
x_{d, n(d)}(t) := \frac{1}{n(d)} \sum_{i=1}^{n(d)} K_d \left(t, \vec{\xi}_i^{(d)} \right) \cdot  f \left(\xi^{(d)}_i \right),  \eqno(1.7)
$$
the so-called Depending Trial Method (DTM), see [9] [10], [20],  and correspondingly consistent as $ \  \min_d n(d) \to \infty \ $ an unbiased Monte-Carlo
estimation  $ \ x_{\vec{n}}^{(M)}(t) $ for $ \ x^{(M)}(t) \ $  may be written as

$$
z_{\vec{n}}^{(M)}(t) := f(t) + \sum_{d=1}^M x_{d, n(d)}(t) =
$$

$$
f(t) + \sum_{d=1}^M \frac{1}{n(d)} \ \sum_{j=1}^{n(d)}  K_d \left(t, \vec{\xi}_i^{(d)} \right) \cdot  f \left(\xi^{(d)}_i \right);  \eqno(1.8)
$$

$ \ \vec{n} := \{ n(1), n(2), \ldots, n(d) \}, \ $ so that as $ \  M \to \infty  \ $ and as $  \  \min_{d} n(d) \to \infty \  $

$$
z_{\vec{n}}^{(M)}(t) \to z^{(M)}(t)
$$
with probability one in the uniform norm:

$$
{\bf P} \left( \ ||  z_{\vec{n}}^{(M)}(t) - z^{(M)}(t) || \to 0 \ \right) = 1,
$$
the Law of Large Numbers (LLN) in the Banach space $ \ C(T). \ $ \par
 \ This method was investigated in [20], see also [10], [4], [5], [8],  [11]-[12], [13], [17], [21] etc. It was proved in particular in [20] that

$$
 \Var_M(\vec{n}) \stackrel{def}{=} \max_t \Var  \left\{ z_{\vec{n}}^{(M)}(t)  \right\} \le ||f||^2 \cdot \sum_{d=1}^M [n(d)]^{-1}  (\rho_2)^d,  \eqno(1.9)
$$
including (formally) the case $ \ M = \infty. \ $ \par

\vspace{4mm}

 \ We introduced  above  the value

$$
\rho_2 = \rho_2[K] \stackrel{def}{=} \max_{t \in T} \int_{T} (K(t,s) )^2 \ \mu(ds), \eqno(1.10)
$$
so that $ \ \rho_2[K] \ $ represents the norm: $ \ \rho_2[K] = ||K^{(2)} ||(C(T) \to C(T)) \ $  of a Kroneker's square $ \  K^{(2)}[\cdot] \  $
of the integral operator $ \ K, \  $, i.e. the compact linear integral operator acting inside the space $ \ C(T) \ $ with the  continuous
kernel $ \  K^2(t,s):  \ $

$$
 K^{(2)}[g] (t) \stackrel{def}{=}  \int_T K^2(t,s) \ g(s) \ \mu(ds).  \eqno(1.11)
$$
 \ Of course,  $ \ \rho \le  [\overline{\rho_2}]^2. $\par

 \ More detail information about integral operators may be found in the classical books  of N.Dunford and J.Schwartz [6], [7]. \par

\vspace{4mm}

 \   One can calculate the general amount $ \ N \ $ of elapsed independent random variables  with  distribution $ \ \mu \ $  as follows

$$
N = \sum_{d=1}^M d \cdot n(d). \eqno(1.12)
$$

\ Solving the following conditional extremal problem

$$
 \Var_M(\vec{n}) \to \min
$$
subject to the restriction

$$
 \sum_{d=1}^M d \cdot n(d) \le N, \eqno(1.13)
$$
the authors of the article [20] obtained the following estimate

$$
\min \Var_M(\vec{n})/ \left[\sum_{d=1}^M d \cdot n(d) \le N \right] \asymp \frac{1}{N}. \eqno(1.14)
$$

 \ On the other words, the speed of convergence $ \  z_{\vec{n}}^{(M)}(t) \to z^{(M)}(t)  \  $   is equal to $ \  1/\sqrt{N}, \  $ alike as in the classical Monte - Carlo method.
 At the same is true also in the uniform norm, see   [20]. \par

 \vspace{4mm}

 \ {\bf We want to improve in this report the described before algorithm, namely, we intend  slightly reduce the number of expended random variables with at the same
 exactness.} \par

\vspace{4mm}

 \ {\bf We will ground also the convergence of our approximation in the uniform norm still with the classical rate and describe the building of confidence region in this norm,
 asymptotical or not. }\par

\vspace{4mm}

 \ Some modern results about the Monte-Carlo solutions of inequations be found in an articles [8], [11]-[12], [13], [20];
see also a monograph of Prem K.Kythe, Pratap Puri [17]. \par

\vspace{4mm}

 {\it Throughout this paper, the letters  $ \ C, C_j(\cdot) \ $ etc. will denote a various positive finite
constants which may differ from one formula to the next even within a single string
of estimates and which does not depend on the essentially variables  $  \ t, x, y, u \ $ etc. \par

\ We make no attempt sometimes to obtain the best values for these constants.}\par

 \ One of the  new applications, namely, in the reliability theory, of these equations  may be found for instance in [10]; see also the
reference therein. \par

 \ Another denotation: we define  for arbitrary random variable $ \ \xi \ $  it {\it centering} $ \ \xi^{(0)} \ $ as an ordinary {\it linear} operator

$$
\xi^{(0)} \stackrel{def}{=} \xi - {\bf E} \xi, \eqno(1.15)
$$
 the ``pure random part'',  which is defined for arbitrary r.v. from the space $ \ L_1(\Omega, {\bf P}), \ $
so that $ \ {\bf E} \ \xi^{(0)} = 0 $ and $ \ \Var \left[ \xi^{(0)} \right] = \Var \ \xi. $ \par

\vspace{4mm}

\section{ Structure of offered solution. Rough estimate of convergence. }

 \vspace{4mm}

 \ Let us introduce the following (functional) recursion:  $ \   x_{0,0}(t)  := f(t), \ t \in T  \  $ and for the values  $ \  m =0,1,2, \ldots, M - 1, \ $ where $ \ M = \const \ge 2 \ $

$$
x_{m+1, n(m+1)}(t) := f(t) + \frac{1}{n(m+1)} \cdot \sum_{l=1}^{n(m+1)} K \left(t, \xi_l^{(m)} \right) \ x_{m, n(m)} \left( \ \xi_l^{(m)} \  \right), \eqno(2.0)
$$
 which is in turn  some modification of the Depending Trial Method. \par

 \vspace{4mm}

 {\it \ We describe here a rough investigations of this approach; the rigorous reasoning will be represented in the next sections. } \par

 \vspace{4mm}

 \ Note first of all that the amount of all elapsed random variables with distribution $ \ \mu \ $ in (2.0), which we will denote again by $ \ N, \ $ is

$$
N = \sum_{d=1}^M n(d), \eqno(2.1)
$$

 \ As before, $ \  \{ \xi_l^{(m)} \} \  $ are common independent random variables with distribution $ \ \mu \ $ and $ \ \{n(m)\} \ $ are certain sequence of natural positive numbers.\par

 \vspace{4mm}

 \ {\it It will be presumed furthermore  that}

 $$
  \forall m= 1,2, \ldots, M \ \Rightarrow n(m) \ \to \infty.  \eqno(2.2)
 $$

\vspace{4mm}

 \ All the functions $ \ x_{m, n(m)}(t), \ m = 1,2,\ldots, M \ $ are  random processes (r.p.), or more generally random fields (r.f). For instance,

$$
 x_{1, n(1)}(t) = f(t) + \frac{1}{n(1)} \sum_{l=1}^{n(1)} K \left( \ t, \xi_l^{(1)} \ \right) \ f \left( \ \xi_l^{(1)} \right). \eqno(2.3)
$$

 \ The last expression may be rewritten under some conditions, which be clarified below, as follows.

$$
 x_{1, n(1)}(t) = x_1(t)  + \frac{1}{\sqrt{n(1)}} \tau_{1, n(1)}(t), \eqno(2.4)
$$
where the r.f. \ $ \tau_{1, n(1)}(t) \ $ is uniformly relative the index $ \ n(1) \ $ subgaussian:

$$
\sup_{n(1) \ge 1} {\bf P} \left( \ \sup_{t \in T} \left| \ \tau_{1, n(1)}(t) \ \right| > u \  \right) \le \exp \left( - C_1 \ u^2 \  \right), \  \eqno(2.5)
$$
 where $ \ C_1 =  C_1[f,K] =\const > 0. \ $ As $ \ n(1) \to \infty \ $ the sequence of r.f. $ \ \{ \tau_{1, n(1)}(t) \} \ $ converges weakly in the space $ \ C(T) \ $
by virtue of CLT in this space   to the
centered continuous a.e. Gaussian random field $ \ \tau_1(t) = \tau_{1, \infty}(t)  \ $   with covariation function

$$
R_1(t_1,t_2)  = \Cov \left( \ \tau_1(t_1), \ \tau_1(t_2)  \ \right) = {\bf E} \tau_1^0(t_1) \ \tau_1^0(t_2) =
$$

$$
\int_T K(t_1,s) \ K(t_2,s) \ f^2(s) \ \mu(ds) -
$$

$$
\int_T K(t_1,s) \ f(s) \ \mu(ds) \cdot \int_T K(t_1,s) \ f(s) \ \mu(ds), \eqno(2.6)
$$
and wherein

$$
\sup_{t \in T} \Var \left[ \ x_{1, n(1)}(t) \ \right] \le \frac{D(1)}{n(1)}, \ D(1) = \const = \max_{t \in T} R_1(t,t) < \infty. \eqno(2.7)
$$

\vspace{4mm}

 \ Let us introduce the following covariation functions, more exactly, the sequence of ones

$$
R_{m+1}(t_1. t_2) \stackrel{def}{=}  \int_T K(t_1,s) \ K(t_2,s) \ x_m^2(s) \ \mu(ds)  -
$$

$$
\int_T  K(t_1,s)  \ x_m(s) \ \mu(ds)  \cdot \int_T  K(t_2,s)  \ x_m(s) \ \mu(ds), \ m = 0,1,2,\ldots. \eqno(2.8)
$$

 \ Define also for arbitrary (continuous) covariation function $  \  R = R(t_1, t_2)  \ $  its linear transform (operator)

$$
V_K[R](t_1,t_2) \stackrel{def}{=} \int_T K(t_1,s) \ K(t_2,s) \ R(s,s) \ \mu(ds). \eqno(2.9)
$$

 \ Evidently, the function  $ \  (t_1, t_2) \to V_K[R](t_1,t_2), \ t_1,t_2 \in T \  $ is also a continuous covariation function. \par

\vspace{4mm}

 \ One can substitute the expression (2.4) for the (random) function  $ \ x_{1, n(1)}(t) \ $ into the recursion (2.0) in order to obtain
the  value for  $ \ x_{2, n(2)}:  \ x_{2, n(2)}(t) =  \ $

$$
 f(t) + \frac{1}{n(2)} \cdot \left\{ \sum_{j=1}^{n(2)} K \left(t, \xi^{(2)}_j \right) \ \left[ \ x_{1}( \xi^{(2)}_j ) + n(1)^{-1/2} \zeta_{1, n(1)} ( \ \xi^{(2)}_j ) \ \right]    \right\} =
$$

$$
 f(t) + \frac{1}{n(2)} \cdot \left\{ \sum_{j=1}^{n(2)} K \left(t, \xi^{(2)}_j \right) \ \ x_{1}( \xi^{(2)}_j ) \right\} +
$$

$$
 n(2)^{-1} \ n(1)^{-1/2}   \left\{ \sum_{j=1}^{n(2)} K \left(t, \xi^{(2)}_j \right) \ \zeta_{1, n(1)} ( \ \xi^{(2)}_j ) \ \right\} =
$$

$$
x_2(t) + \frac{1}{ \sqrt{n(2)}} \tau_2(t) + \frac{1}{ \sqrt{ n(2) \ n(1)}} \tau_{2,1}(t),
$$
where the centered non - correlated r.f. $ \ \tau_2(t), \  \tau_{2,1}(t)  \  $  have correspondingly covariation functions

$$
\Cov( \tau_2(t_1), \  \tau_2(t_2) ) = R_2(t_1,t_2)
$$
and

$$
R_{2,1}(t_1,t_2) \stackrel{def}{=}
\Cov( \tau_{2,1}(t_1), \  \tau_{2,1}(t_2) ) = V_K[R_1](t_1,t_2). \eqno(2.10)
$$

 \ Let us denote for simplicity

$$
\sigma(m+1) := \frac{1}{\sqrt{n(m+1)}}, \ \sigma(m+1,m) := \frac{1}{ \sqrt{n(m+1) \ n(m)}}, \
$$

$$
\sigma(m+1,m,m-1) := \frac{1}{ \sqrt{n(m+1) \ n(m) \ n(m-1)}}, \ \ldots,
$$

$$
\sigma(m+1,m,m-1, \ldots,m-k) =
$$

$$
\frac{1}{ \sqrt{n(m+1) \ n(m) \ n(m-1) \ n(m-2) \ \ldots \ n(m-k)}}, \ k = 0,1,\ldots,m-1;
$$
so that

$$
\sigma(m+1,m,m-1, \ldots,1) = \frac{1}{ \sqrt{ \prod_{l=-1}^{m-1} n(m-l)}}.
$$

\ We  deduce  passing to the more general case the following decomposition:

$$
x_{m+1, n(m+1)}(t) = x_{m+1}(t) +
$$

$$
 \sigma(m+1) \ \tau_{m+1}(t)  + \sigma(m+1,m) \ \tau_{m+1,m}(t)  +
$$

$$
\sigma(m+1,m,m-1) \ \tau_{m+1, m, m-1}(t) + \ldots  +
$$

$$
\sigma(m+1,m,m-1, \ldots,1) \ \tau_{m+1, m, m-1, \ldots, 2,1}(t), \eqno(2.11)
$$
where  the  centered and non - correlated r.f. $ \ \{\tau_{\vec{k}}(t) \ \} \ $ have the covariation function correspondingly

$$
\Cov \left( \tau_{m+1}(t_1), \  \tau_{m+1}(t_2) \  \right) = R_{m+1}(t_1, t_2);
$$

$$
R_{m+1,m} (t_1,t_2):= \Cov \left( \tau_{m+1,m}(t_1), \  \tau_{m+1,m}(t_2) \  \right) = V_K \left[R_{m,m-1} \right](t_1,t_2);
$$

$$
R_{m+1,m,m-1}(t_1,t_2) :=
 \Cov \left( \tau_{m+1,m, m-1}(t_1), \  \tau_{m+1,m, m-1}(t_2) \  \right) = V_K \left[R_{m,m-1,m-2} \right](t_1,t_2);
$$

\vspace{3mm}

.........................................................................................................................\\

$$
R_{m+1,m,m-1,\ldots,1}(t_1,t_2):=
 \Cov \left( \tau_{m+1,m, m-1,\ldots,1}(t_1), \  \tau_{m+1,m, m-1,\ldots,1}(t_2) \  \right) =
$$

$$
  V_K \left[R_{m,m-1,m-2,\ldots,1} \right](t_1,t_2), \ - \eqno(2.12)
$$
 a recursion: $ \ m = 2,3,\ldots, M-1. \ $ \par

\vspace{3mm}

 \ Of course, all the  centered continuous random fields $ \ \tau_{m+1, m, m-1, \ldots, m-k} (t), \ 1 \le k \le m-1, \ m = 1,2, \ldots,M  \ $ dependent on the vector $ \ \vec{n}_m = \vec{n} = \\
 \{ n(m), n(m-1), \ldots, n(1)  \}: \ $

$$
\tau_{m+1, m, m-1, \ldots, m-k} (t) = \tau^{(\vec{n})}_{m+1, m, m-1, \ldots, m-k} (t).
$$

\  They are uniformly subgaussian:

$$
\sup_{\vec{n}} \  {\bf P} \left( \sup_{t \in T} \left| \ \tau_{m+1, m, m-1, \ldots, m-k} (t) \ \right| > u  \right) \le
$$

$$
\exp \left( -C(K(\cdot,\cdot), f(\cdot); \ m+1, m, m-1, \ldots, m-k) \ u^2 \  \right), \eqno(2.13)
$$
$ \  C(K(\cdot,\cdot), f(\cdot); \ m+1, m, m-1, \ldots, m-k) > 0, \ $
 and converges weakly in distribution
 in the space $ \ C(T) \ $ as $ \ \min_m n(m) \to \infty \ $ to the continuous centered independent Gaussian random fields with covariations
 described in (2.12). \par

\ Note that the variance $ \  \Var \left[ \ x_{M, n(M) }(t) \ \right]   \ $ allows the following estimate

$$
C^{-1}(K,f) \ \max_{t \in T} \Var \left[ \ x_{M, n(M) }(t) \ \right] \le \frac{1}{n(M)} + \frac{1}{n(M) \ n(M-1)} +
$$

$$
\frac{1}{n(M) \ n(M-1) \ n(M-2) } + \ldots + \frac{1}{n(M) \ n(M-1) \ n(M-2) \ \ldots \ n(2) \ n(1)}, \eqno(2.14)
$$
$ \ C(K,f) \in (0,\infty), \  $ and besides if for instance both the functions $ \  f(\cdot), \ K(\cdot, \cdot) \  $ are positive, the inverse inequality  to the one in (2.14) holds true
relative the other constant. \par

\vspace{3mm}

 \ One can choose for instance $ \ n(M):= N/2, \ $ and further

$$
n(M-1) := N/4, \ n(M-2) := N/8, \ \ldots, n(M-k) := N/2^{k+1}, \eqno(2.15)
$$
 to make sure  that the offered algorithm in  comparison  to the previous described before obeys  at the same speed of
 convergence  but requires significantly less random variables, cf. the relations (2.1) and (1.12). \par
 \ {\it It will be presumed in the sequel of course that }

 $$
   N >> 2^{M+1}.   \eqno(2.16)
 $$

\vspace{4mm}

\section{Central Limit Theorem for described solution. \\ Asymptotical confidence interval in the unifirm norm.}

 \vspace{4mm}

 \ We recall now some notions from the theory of the Central Limit Theorem (CLT) in the separable Banach space of continuous functions $ \ C(T). \ $
 More detail information about this theory  may be found in [2], [5], [15], [18], [22] etc. \par

  \ Let $ \ \nu = \nu(t), \ t \in T \ $ be centered: $ \ {\bf E} \nu(t) = 0, \ t \in T \ $  numerical valued random field having finite second moment:

$$
\forall t,s \in T \ \exists R(t,s)  = \Cov(\nu(t), \nu(s)). \eqno(3.0)
$$
 \ Let also $ \ \nu_i(t) \ $ be independent copies of $ \ \nu(t). \ $ Denote

$$
S_n(t) := n^{-1/2} \sum_{i=1}^n \nu_i(t). \eqno(3.1)
$$
 \  The r.f. $ \ \nu(\cdot) \ $ satisfies by definition CLT in the space $ \ C(T), \ $ if r.f. $ \ \nu(t) \ $ is continuous a.e. and the sequence of r.f. $ \ S_n(\cdot) \ $
converges as $ \ n \to \infty \ $ weakly in  distribution to the Gaussian r.f. $ \  S_{\infty}(t) = S(t),  \ $ also continuous everywhere, having at the same first two moments as the source
r.f. $ \nu(t). \ $ \par

 \ This implies that for arbitrary continuous bounded functional $ \ F: C(T) \to R \ $

$$
\lim_{n \to \infty} {\bf E} F(S_n) = {\bf E} F(S). \eqno(3.2)
$$

 \ As a consequence from (3.2)

$$
\lim_{n \to \infty} {\bf P}(\max_{t \in T} |S_n(t)| > u) = {\bf P}(\max_{t \in T} |S(t)| > u), \ u > 0. \eqno(3.3)
$$

 \ Note that the asymptotical as $ \ u \to \infty \ $  behavior as well as non-asymptotical at $ \ u \ge 1  \ $  estimates  are well known, see e.g. [21], [18].
 This circumstance  allows to construct a confidence region in the Monte-Carlo parametric computations, see [9], [10], [18]. \par

 \vspace{4mm}

 \ We must recall also some used facts about  the so - called subgaussian random variables. The r.v. $ \ \xi \ $ defined on some
probability space is said to  be subgaussian, write: $ \ \xi \in \Sub, \ $ if there exists a non-negative constant $ \ \tau = \tau(\xi) =
\tau(\Law(\xi)), \ $ for which

$$
\forall \lambda \in R \ \Rightarrow \max_{\pm} {\bf E} \exp( \pm \lambda \xi) \le \exp(0.5 \ \lambda^2 \tau^2). \eqno(3.4)
$$

 \ This notion was introduced at first by Kahane J.P. in [14].  See also [2], [3], [5], [10], [18], [19], [22] etc. \par

 \  Evidently, if $ \ \xi \in \Sub, \ \xi \ne 0,  \ $ then $ \ {\bf E} \xi = 0 \ $ and

$$
 \max \left[ {\bf P}(\xi \ge u),  {\bf P}(\xi \le -u)    \right] \le
\exp \left( - 0.5 \ u^2 /\tau^2 \right), \ u \ge 0; \eqno(3.5)
$$
and inverse conclusion up to multiplicative constant is also true.\par

 \ The minimal non-negative value $ \ \tau = \tau(\xi) \ $  satisfying (3.4) for all the values $ \ \lambda \in R, \ $  is named as
{\it subgaussian norm \ }  of the r.v. $ \ \xi, \ $  write $ \ ||\xi||\Sub: \ $

$$
 ||\xi||\Sub \stackrel{def}{=}
\max_{\pm} \sup_{\lambda \ne 0} \left\{ \left[ \ln {\bf E} \exp(  \ \pm \ \lambda \ \xi)  \   \right]^{1/2}/|\lambda| \ \right\}. \eqno(3.6)
$$

 \   Buldygin V.V. and  Kozachenko Yu.V. in [2], see also [3] proved
that the functional $ \  \xi \to ||\xi||\Sub \  $ is really the norm and the space $ \ \Sub\ $  forms the complete rearrangement invariant (r.i.)
in the  classical sense [1], chapters 1,2 Banach space. For  instance, if $ \  \xi_1, \xi_2 \in \Sub, \ c_1, c_2 = \const \in R, \ $ then

$$
||c_1 \xi_1 + c_2 \eta||\Sub  \le |c_1 | ||\xi||\Sub + |c_2| \ ||\eta|| \Sub, \ c_1, c_2 = \const \in R.
$$
 \ If in addition both the r.v. $ \ \xi, \ \eta \ $ are independent, then

$$
||c_1 \xi_1 + c_2 \eta||\Sub  \le  \left[ |c_1^2 | ||\xi||^2\Sub + |c_2|^2 \ ||\eta||^2 \Sub \right]^{1/2}, \eqno(3.7)
$$
and  analogously for the linear combination of several independent random variables.  In particular, if $ \ \eta_i, \ i = 1,2,\ldots,n; \ \eta = \eta_1 \ $
are identical distributed centered subgaussian random variables, then

$$
\sup_n ||n^{-1/2} \ \sum_{i=1}^n \eta_i|| \Sub = ||\eta||\Sub. \eqno(3.7a)
$$

 \ The centered Gaussian distributed r.v. $ \ \xi \ $  is also subgaussian and wherein $ \ ||\xi||\Sub = [\Var(\xi)]^{1/2}. \ $ A more interest
fact: let  the {\it mean zero} r.v. $ \ \eta \ $  be bounded: $ \ ||\eta||_{\infty} < \infty; \ $ then it is also subgaussian and herewith
	
$$
||\eta||\Sub \le ||\eta||_{\infty}, \eqno(3.8)
$$
and the last estimate is in general case non-improvable. \par
 \ For instance, the Rademacher distributed r.v. $ \ \eta: \ {\bf  P}(\eta = 1) = {\bf  P}(\eta = -1)  = 1/2 \ $ is  subgaussian and
$ \ ||\eta||\Sub = 1. \ $ \par

\vspace{4mm}

  \ Let us return to the source r.f. $ \ \nu = \nu(t), \ t \in T. \ $ Suppose now that it is uniformly subgaussian relative the parameter $ \ t: \ $

$$
D = D_{\nu} := \sup_{t \in T} ||\nu(t)||\Sub < \infty.
$$

 \ We introduce so-called {\it  natural} finite semi-distance function $ \ \Delta(t,s) = \Delta_{\nu}(t,s) \ $ as follows

$$
 \Delta(t,s) = \Delta_{\nu}(t,s)  \stackrel{def}{=} ||\nu(t) - \nu(s)||\Sub. \eqno(3.9)
$$
 \ Evidently,

$$
\diam_{\Delta}T \stackrel{def}{=}  \sup_{t,s \in T} \Delta(t,s) \le 2 D; \eqno(3.10)
$$
and the correspondent  metric entropy   function $ \ H(T, \Delta, \epsilon)  \ $ for the set   $ \ T \ $ relative the distance  $ \ \Delta[\nu](t_1, t_2)  \ $
at the point $ \ \epsilon, \ \epsilon \in (0, \diam_{\Delta}T). $ \par

\vspace{4mm}

 \ {\bf Lemma 3.0; \ see  e.g. [5], [15], [18] etc.}  If for the subgaussian random field

$$
 \nu = \nu(t), \ t \in T  \ \Rightarrow \  D_{\nu} < \infty \ \eqno(3.11a)
$$

and  the following  so-called {\it entropic} integral converges:

$$
 I = I(\nu) :=   \int_0^{\diam_{\Delta}T} H^{1/2}(T, \Delta_{\nu}, \epsilon) \ d \epsilon < \infty, \eqno(3.11b)
$$
then the r.f. $ \ \nu(t) \ $  satisfies CLT in the space $ \ C(T) \ $ and moreover

$$
\sup_n {\bf P} \left( \ \sup_{t \in T} | S_n(t)| > u \  \right) \le \exp \left(- \ C_2(I) \ (u/D)^2 \  \right), C_2 > 0. \eqno(3.12)
$$

\vspace{4mm}

 \ {\bf Example 3.1.}  \ Note that the condition (3.11b) is fulfilled if for example $ \ T \ $ is bounded subset of the whole Euclidean space $ \ R^d \ $ equipped
with ordinary distance $ \ ||t - s||  $ and if

$$
\Delta_{\nu}(t,s) \le C \ ||t-s||^{\alpha}, \ \exists \alpha \in (0,1], \ \exists  C \in (0,\infty), \ \eqno(3.12a),
$$
H\"older's condition. \par

 \ In order to formulate and prove the CLT for our solution, we must introduce some new notations.

$$
\beta = \beta(K) := \max_{t_1, t_2 \in T} \int_T |K(t_1,s) \ K(t_2,s) | \ \mu(ds) < \infty, \eqno(3.13)
$$

$$
d = d(t_1, t_2) = d^{(K)}(t_1,t_2) := \left\{ \int_T \ \left[ \ K(t_1, s) - K(t_2,s) \ \right]^2 \ \mu(ds) \  \right\}^{1/2}, \eqno(3.14)
$$

$$
d_{m+1}(t_1, t_2) = d_{m+1}^{(K)}(t_1, t_2) :=
$$

$$
\left\{ \int_T \ \left[ \ K(t_1, s) - K(t_2,s) \ \right]^2 \  x_m^2(s) \ \mu(ds) \ \right\}^{1/2}. \eqno(3.14a)
$$

\vspace{4mm}

 \ All the introduced functions  $ \  d(t_1, t_2), \  d_l(t_1, t_2) \ $ are finite semi-distance functions defined on the set $ \ T^2. \ $

\vspace{4mm}

 \ {\bf Lemma 3.1.}  Suppose  that

$$
\int_0^1 H^{1/2}(T, d, \epsilon) \ d \epsilon < \infty. \eqno(3.15)
$$
 \ Then each  centered random fields mentioned in (2.11):  $ \  \tau_{m+1}(t), \ \tau_{m+1,m}(t),    \ \tau_{m+1, m, m-1}(t), \ldots,  \ \tau_{m+1, m, m-1, \ldots, 2,1}(t) \ $
satisfies the CLT as $ \min_m n(m) \to \infty \ $ in the space $ \ C(T) \ $ with covariation functions described correspondingly in (2.12). \par

\vspace{4mm}

 \ {\bf Proof.} We deduce

$$
|R_1(t,s)| \le \beta \ ||f||^2, \ R_2(t,s)| \le \beta \ ||x_1||^2
$$
and in general

$$
 |R_{m+1}(t,s)| \le \beta \ ||x_m||^2, \ m = 1,2,\ldots.
$$
 \ As long as

$$
||x_m|| \le \frac{||f||}{1-\rho}, \ m \ge 1,
$$
we obtain

$$
|R_{m+1}(t,s)| \le \beta \ \frac{||f||^2}{(1 - \rho)^2}. \eqno(3.16)
$$

 \ Further,

$$
|R_{m+1,m}(t,s)| \le \beta \ \max_{t,s} |R_{m,m-1}(t,s)|,
$$
therefore $ \ R_{(M)} \stackrel{def}{=} $

$$
 \max_{t,s \in T} \left\{ \ |R_{M+1}(t,s)|  + |R_{M+1,M}(t,s)| +  \ \ldots  +
 |R_{M+1,M,M-1, \ \ldots, 1}(t,s)| \ \right\} \le
$$

$$
 ||f||^2 \ \frac{\beta + \beta^2 + \ldots \ + \beta^M}{(1 - \rho)^2} =
||f||^2 \ \frac{|\beta - \beta^{M+1}|}{|1-\beta| \ (1 - \rho)^2}. \eqno(3.17)
$$
 \ Evidently, if $ \ \beta < 1, \ $ then

$$
R_{(M)} \le ||f||^2 \ \frac{\beta}{(1 - \beta) \ (1 - \rho)^2}. \eqno(3.17a)
$$

 \ We continue to estimate:

$$
d_{m+1}(t_1, t_2) \le  \frac{||f| \ d(t_1,t_2)}{1 - \rho},
$$
and analogously $ \tilde{d}_{m,k}(t_1,t_2) :=  $

$$
d_{R(m+1,m, \ldots,m-k)}(t_1,t_2) \stackrel{def}{=} || \tau_{m+1,m,m-1,\ldots,m-k}(t_1) - \tau_{m+1,m,m-1,\ldots,m-k}(t_2)||\Sub \le
$$

$$
\frac{||f|| \ d(t_1,t_2)}{(1-\rho)^k}. \eqno(3.18)
$$

 \ On the other words, both the distance functions  $ \ d_{R(m+1,m, \ldots,m-k)}(t_1,t_2) \ $ and $ \ d(t_1, t_2) \ $ are linear equivalent. Following, the
 entropic integral for the metric  $ \ d_{R(m+1,m, \ldots,m-k)}(t_1,t_2) \ $  is  finite as well as one  for the metric $ \ d, \ $ see (3.15).\par
  \ The finiteness of the diameter of the set $ \ T \ $ relative both the considered distances follows straightforwardly  from the estimate (3.17). \par
 \ It remains to apply the proposition of Lemma 3.0. \par

\vspace{3mm}

 \ Let us formulate and prove one of the main results of this report. \par

\vspace{4mm}

 \ {\bf Theorem 3.1. CLT in the space of continuous functions for our solutions.} \par

\vspace{4mm}

 \ Assume that all the formulated before conditions, including the restrictions (2.15) and (2.16) are fulfilled. We propose that the sequence of random fields

$$
\gamma_n (t)= \gamma_{\vec{n}}(t) := \sqrt{n(M)} \left( \ x_{M,n(M)}(t) - x_M(t) \ \right) \eqno(3.19)
$$
converges weakly in distribution as $ \ n(M) \to \infty \ $ in the space $ \ C(T) \ $ to the continuous a.e centered Gaussian random field $ \ \gamma(t) \ $ with
covariation function

$$
\Cov(\gamma(t), \ \gamma(s)) = R_{M}(t,s). \eqno(3.20)
$$

\vspace{4mm}

 \ {\bf Proof.} We exploit the decomposition (2.11):  $ \ x_{M, n(M)}(t) = x_{M}(t) + $

$$
 \sigma(M) \ \tau_{M}(t)  + \sigma(M,M-1) \ \tau_{M,M-1}(t)  + \sigma(M,M-1,M-2) \ \tau_{M,M-1,M-2}(t) + \ldots =
$$

$$
 x_M(t) +\sigma(M) \ \tau_{M}(t) + \Sigma_2 = \Sigma_1 + \Sigma_2. \eqno(3.21)
$$
 \ With regard to the first member $ \ \Sigma_1: \ $  it is centered  satisfied the CLT in the space  $ \ C(T) \ $  with parameters showed in (3.20);
we need to ground that the second  member $ \ \Sigma_2 \ $ in (3.21) tends to zero in the uniform  norm with probability one. \par
 \ It is sufficient to consider only the next member in (2.11); the remains ones may be investigated analogously. \par

\ We observe using (2.15), (2.16) and Lemma 3.1

$$
P_N(u) := {\bf P} ( ||x_{M, M-1}(\cdot)|| > u ) = {\bf P} \left(  (n(M) \ n(M-1))^{-1/2} \ ||\tau_{M,M-1}(\cdot) || > u  \right) \le
$$

$$
\exp \left( \  - \ C_4 \ N \ u^2 \  \right), \ u > 0. \eqno(3.22)
$$
 \ Observe that

$$
\forall u > 0 \ \Rightarrow \ \sum_{N > 1} P_N(u) < \infty,
$$
 therefore

$$
{\bf P} \left( || x_{M, M-1}(\cdot)|| \to 0 \right) = 1,
$$
by virtue of Lemma Borel-Cantelli. \par

\vspace{4mm}

 \ Note that if the numbers $ \ n(K) \ $ are choosed in accordance with relation (2.16), then the sequence od r.f.

$$
\gamma_n (t)= \gamma_{\vec{n}}(t) := \sqrt{N/2} \cdot \left( \ x_{M,n(M)}(t) - x_M(t) \ \right) \eqno(3.23)
$$
converges weakly in distribution as $ \ N \to \infty \ $ in the space $ \ C(T) \ $ to the continuous a.e centered Gaussian random field $ \ \gamma(t) \ $ with at the same
covariation function

$$
\Cov(\gamma(t), \ \gamma(s)) = R_{M}(t,s). \eqno(3.24)
$$

\vspace{4mm}

\section{ Non-asymptotic confidence region in the uniform norm.}

\vspace{4mm}

 \  We construct in this section the non-asymptotic confidence interval  for the $ \ x_{M} (t). \ $ \par

\vspace{4mm}

 \ {\bf Theorem 4.1.}  The following non-asymptotical estimate is valid under at the same assumptions as in the foregoing theorem 3.1:

$$
\sup_{N \ge 4}  {\bf P} \left( \ \sqrt{N/2} \cdot  \sup_{t \in T} \left| \ x_{M,n(M)}(t) - x_M(t) \ \right| > u \right) \le
$$

$$
\exp(-C_3(K,f;M)\ u^2), \ u > 0. \eqno(4.1)
$$

\vspace{4mm}

{\bf Proof.} Let us consider the following r.f., more precisely, the sequence of random fields

$$
\upsilon_N(t) = \upsilon_{N,M, \vec{n}} (t) := \sqrt{n(M)} \left(  \ x_{M,n(M)}(t) - x_M(t)  \right), \ t \in T.
$$
 \ It follows immediately from the estimates  (3.17), (3.18) that

$$
\sup_{N \ge 4} \ \sup_{t \in T} \ ||\upsilon_N(t)||\Sub = C_4 < \infty, \eqno(4.2)
$$

$$
\sup_{N \ge 4} \  ||\upsilon_N(t) -  \upsilon_N(s)||\Sub = C_5 \ d(t,s) < \infty. \eqno(4.3)
$$

 \ Our proposition (4.1) follows from one in the Lemma 3.0. \par

\vspace{4mm}

\section{ Concluding remarks.}

\vspace{4mm}

 \ {\bf A.} The offered in this preprint method may be easily generalized on the systems of integral equations, and perhaps
on the some non-linear ones. \par

\vspace{4mm}

 \ {\bf B.} Perhaps, the offered here method  may be used for the solving by the Monte-Carlo method for integral equations containing
kernels  discontinuously depending on some parameter (parameters), in the spirit of the article [10].  \par

 \vspace{4mm}

\begin{center}

\ {\bf References.}

\end{center}

 \vspace{4mm}

{\bf 1. Bennet C., Sharpley R.}  {\it  Interpolation of operators.} Orlando, Academic
Press Inc., (1988). \\

 \vspace{4mm}

{\bf 2.  Buldygin V.V., Kozachenko Yu.V. }  {\it Metric Characterization of Random
Variables and Random Processes.} 1998, Translations of Mathematics Monograph, AMS, v.188. \\

\vspace{4mm}

{\bf 3.  Buldygin V.V., Kozachenko Yu.V. }  {\it About subgaussian random variables.} Ukrainian Math. Journal, 1980, V. 32, Issue 6,
pp. 723 - 730. \\

 \vspace{4mm}

{\bf 4.  Arnaud Doucet, Adam M. Johansen, and  Vladislav B. Tadic.} {\it On solving integral equations using Markov chain Monte Carlo methods.}
Applied Mathematics and Computation 216 (2010) 2869–2880. \\

\vspace{4mm}

{\bf 5. Dudley R.M.} {\it Uniform Central Limit Theorem.} Cambridge, University Press, (1999), 352-367.\\

 \vspace{4mm}

{\bf 6.  Dunford N., Schwartz J.}  {\it Linear operators. Part 1: General Theory.} (1958),
Interscience Publishers, New York, London. \\

\vspace{4mm}

{\bf 7. Dunford N., Schwartz J.}  {\it Linear operators. Part 2: Spectral Theory.} (1963),
Interscience Publishers, New York, London.\\

 \vspace{4mm}

 {\bf 8. R.Farnoosh, M.Ebrahimi.} {\it Monte Carlo method for solving Fredholm integral equations of the second kind.}
Applied Mathematics and Computation, Volume 195, Issue 1, 15 January 2008, Pages 309-315. \\

\vspace{4mm}

{\bf  9. Frolov  A.S.,  Chentzov  N.N.}
{\it On the calculation by the Monte-Carlo method of definite integrals depending on the parameters.}
Journal of Computational Mathematics and Mathematical Physics, (1962), V. 2, Issue 4, p. 714-718, (in Russian).\\

\vspace{4mm}

{\bf 10.  Grigorjeva M.L., Ostrovsky E.I.} {\it Calculation of Integrals on discontinuous functions by means of depending trials
method.} Journal of Computational Mathematics and Mathematical Physics, (1996), V. 36, Issue 12,
p. 28-39 (in Russian).\\

\vspace{4mm}

 {\bf 11. Hong Zhi Min , Yan Zai Zai and Chen Jian Rui.} {\it  Monte Carlo Method for Solving the Fredholm Integral Equations of the Second Kind.}
  Transport Theory and Statistical Physics, Volume 41, 2012-Issue 7, pp.144 - 153.\\

 \vspace{4mm}

{\bf 12. Zhimin Hong, Xiangzhong Fang, Zaizai Yan1, Hui Hao.  }
{\it On Solving a System of Volterra Integral Equations with Relaxed Monte Carlo Method.}
JAMP, Vol.4 No.7, July 2016 pp.1315 - 1320.\\

\vspace{4mm}

{\bf 13. A. M. Johansen.}  {\it Monte Carlo Solution of Integral Equations (of the second kind).} March 7th, 2011. Internet publication, PDF. \\

  \vspace{4mm}

{\bf 14. Kahane J.P.} {\it  Properties locales des fonctions a series de Fourier  aleatoires.  } Studia Math., 1960, B.19, $ N^o \ $ 1, pp. 1-25. \\

\vspace{4mm}

{\bf 15. Kozachenko Yu. V., Ostrovsky E.I. }  (1985). {\it The Banach Spaces of random Variables of subgaussian Type. } Theory of Probab.
and Math. Stat. (in Russian). Kiev, KSU, 32, 43-57. \\

\vspace{4mm}

{\bf 16. Kozachenko Yu.V., Ostrovsky E., Sirota L.}  {\it Relations between exponential tails, moments and
moment generating functions for random variables and vectors.} \\
arXiv:1701.01901v1 [math.FA] 8 Jan 2017 \\

\vspace{4mm}

 {\bf 17. Prem K.Kythe, Pratap Puri.} {\it  Computational Methods for Linear Integral Equations.}
 Springer science + Business Media, LLC, USA, 2002. \\

\vspace{4mm}

{\bf 18. Ostrovsky E.I. } (1999). {\it Exponential estimations for Random Fields and its
applications,} (in Russian). Moscow-Obninsk, OINPE. \\

 \vspace{4mm}

{\bf 19. Ostrovsky E. and Sirota L.} {\it Vector rearrangement invariant Banach spaces
of random variables with exponential decreasing tails of distributions.} \\
 arXiv:1510.04182v1 [math.PR] 14 Oct 2015 \\

 \vspace{4mm}

 {\bf 20. Ostrovsky E. and Sirota L. }  {\it   Monte-Carlo method for multiple parametric integrals
calculation and solving of linear integral Fredholm equations of a second kind, with confidence regions in the uniform norm.} \\
arXiv:1101.5381v1 [math.FA] 27 Jan 2011 \\

\vspace{4mm}

{\bf 21. V.I.Piterbarg.} {\it Asymptotical methods in the theory of Gaussian random processes and fields.} Moscow, 1976, MSU, (in Russian). \\

\vspace{4mm}

 {\bf 22.  O.I.Vasalik, Yu.V.Kozachenko, R.E.Yamnenko.} {\it $ \ \phi \ - $ subgaussian  random processes. } Monograph, Kiev, KSU,
2008;  (in Ukrainian). \\

\vspace{4mm}

\end{document}